\documentclass[12pt]{article}

\usepackage{amsmath, epsfig, cite}
\usepackage{amssymb}
\usepackage{amsfonts}
\usepackage{latexsym}
\usepackage{float}

\newtheorem{thm}{Theorem}[section]

\newtheorem{cor}[thm]{Corollary}
\newtheorem{conj}[thm]{Conjecture}
\newtheorem{lem}[thm]{Lemma}

\newcommand*\pFqskip{8mu}
\catcode`,\active
\newcommand*\pFq{\begingroup
        \catcode`\,\active
        \def ,{\mskip\pFqskip\relax}%
        \dopFq
}
\catcode`\,12
\def\dopFq#1#2#3#4#5{%
        {}_{#1}F_{#2}\biggl[\genfrac..{0pt}{}{#3}{#4};#5\biggr]%
        \endgroup
}

\numberwithin{equation}{section}

\newcommand{\qed}{{\hfill$\square$}\medskip}

\begin{document}


\begin{center}
{\Large\bf Congruences for truncated hypergeometric series ${}_2F_1$}
\end{center}

\vskip 2mm \centerline{Ji-Cai Liu}
\begin{center}
{\footnotesize Department of Mathematics, East China Normal University, Shanghai 200241, PR China\\
{\tt jc2051@163.com} }
\end{center}


\vskip 0.7cm \noindent{\bf Abstract.}
Rodriguez-Villegas conjectured four supercongruences associated to certain elliptic curves, which were first confirmed by Mortenson by using the Gross-Koblitz formula. In this paper, we aim to prove four supercongruences between two truncated hypergeometric series ${}_2F_1$. This generalizes these four Rodriguez-Villegas supercongruences.

\vskip 2mm
\noindent{\it MR Subject Classifications}: Primary 11A07; Secondary 33C05

\section{Introduction}
In 2003, Rodriguez-Villegas \cite{RV} studied hypergeometric families of Calabi-Yau
manifolds. He observed numerically some remarkable supercongruences between the values of the truncated hypergeometric series and expressions derived from the number of $F_p$-points of the associated Calabi-Yau manifolds. A number of supercongruences for hypergeometric Calabi-Yau manifolds have been conjectured by Rodriguez-Villegas. For manifolds of dimension $d=1$, he conjectured four interesting supercongruences associated to certain elliptic curves. These four supercongruences were first confirmed by Mortenson \cite{Mort1,Mort2} by using the Gross-Koblitz formula.

We first define the truncated hypergeometric series.
For complex numbers $a_i, b_j$ and $z$, with none of the $b_j$ being negative integers or zero,
the truncated hypergeometric series are given by
\begin{align*}
\pFq{r}{s}{a_1,a_2,\cdots,a_r}{b_1,b_2,\cdots,b_s}{z}_n=\sum_{k=0}^{n-1}\frac{(a_1)_k (a_2)_k\cdots (a_r)_k}{(b_1)_k (b_2)_k\cdots (b_s)_k}\cdot \frac{z^k}{k!},
\end{align*}
where
$(a)_0=1$ and $(a)_k=a(a+1)\cdots (a+k-1)$ for $k\ge 1$.

\begin{thm}\label{thm1} (Rodriguez-Villegas-Mortenson)
Let $p\ge 5$ be a prime. Then
\begin{align*}
\pFq{2}{1}{\frac{1}{2},\frac{1}{2}}{1}{1}_p\equiv \left(\frac{-1}{p}\right) \pmod{p^2},\quad
\pFq{2}{1}{\frac{1}{3},\frac{2}{3}}{1}{1}_p\equiv \left(\frac{-3}{p}\right) \pmod{p^2},\\
\pFq{2}{1}{\frac{1}{4},\frac{3}{4}}{1}{1}_p\equiv \left(\frac{-2}{p}\right) \pmod{p^2},\quad
\pFq{2}{1}{\frac{1}{6},\frac{5}{6}}{1}{1}_p\equiv \left(\frac{-1}{p}\right) \pmod{p^2},
\end{align*}
where $\left(\frac{\cdot}{p}\right)$ denotes the Legendre symbol.
\end{thm}

For more proofs of Theorem \ref{thm1}, see \cite{CLZ,Szh1,Taur1,Taur2}. Some extensions of Theorem \ref{thm1} to the case modulus $p^3$ were obtained in \cite{Szh2,Szw}. For some interesting $q$-analogues of Theorem \ref{thm1}, one refers to \cite{GPZ,GZ}.
By studying the generalized Legendre polynomials, Z.-H. Sun \cite{Szh1} extended Theorem \ref{thm1} as follows.
\begin{thm}\label{thmsun} (Z.-H. Sun)
Let $p\ge 5$ be a prime. For any $p$-adic integer $x$, we have
\begin{align}
\pFq{2}{1}{x,1-x}{1}{1}_p\equiv (-1)^{\langle -x\rangle_p} \pmod{p^2},\label{sun-1}
\end{align}
where $\langle a\rangle_p$ denotes the least non-negative integer $r$ with $a\equiv r \pmod{p}$.
\end{thm}

Observe that
\begin{align*}
(-1)^{\langle -1/2\rangle_p}=\left(\frac{-1}{p}\right),\quad
(-1)^{\langle -1/3\rangle_p}=\left(\frac{-3}{p}\right),\\
(-1)^{\langle -1/4\rangle_p}=\left(\frac{-2}{p}\right),\quad
(-1)^{\langle -1/6\rangle_p}=\left(\frac{-1}{p}\right).
\end{align*}
Then Theorem \ref{thmsun} reduces to Theorem \ref{thm1} when $x=\frac{1}{2},\frac{1}{3},\frac{1}{4},\frac{1}{6}$.

Ap\'ery introduced the numbers $A_n=\sum_{k=0}^n{n\choose k}^2{n+k\choose k}^2$ in his ingenious proof \cite{Apery} of the irrationality of $\zeta(3)$. These numbers are known as Ap\'ery numbers.  Since the appearance of these numbers, some interesting arithmetic properties have been gradually discovered. For example, Gessel \cite{Gess} proved that for any prime $p\ge 5$
\begin{align*}
A_{np}\equiv A_n\pmod{p^3},
\end{align*}
which confirmed a conjecture by Chowla et al. \cite{CCC}.

In this paper, we aim to prove a similar type of supercongruences for the truncated hypergeometric series ${}_2F_1$, which generalizes these four supercongruences in Theorem \ref{thm1} as follows.

\begin{thm}\label{thm2}
Suppose $p\ge 5$ is a prime and $n$ is a positive integer. For $x\in \{ \frac{1}{2},\frac{1}{3},\frac{1}{4},\frac{1}{6}\}$, we have
\begin{align}
\pFq{2}{1}{x,1-x}{1}{1}_{np}&\equiv (-1)^{\langle -x\rangle_p}\cdot\pFq{2}{1}{x,1-x}{1}{1}_n \pmod{p^2}.\label{rv}
\end{align}
\end{thm}

Theorem \ref{thm2} reduces to Theorem \ref{thm1} when $n=1$.
Replacing $n$ by $p^{r-1}$ in \eqref{rv} and then using induction, we immediately get the following results:
\begin{cor}
Suppose $p\ge 5$ is a prime and $r$ is a positive integer. For $x\in\{ \frac{1}{2},\frac{1}{3},\frac{1}{4},\frac{1}{6}\}$, we have
\begin{align*}
\pFq{2}{1}{x,1-x}{1}{1}_{p^r}\equiv (-1)^{\langle -x\rangle_p\cdot r} \pmod{p^2}.
\end{align*}
\end{cor}

The organization of this paper is as follows.
In the next section, we first recall some properties of the  Fermat quotients and some combinatorial identities involving harmonic numbers, and then prove two congruences. In Section 3, we will give a new proof of Theorem \ref{thmsun} by using combinatorial identities.
The proof of Theorem \ref{thm2} will be given in the Section 4. We make some concluding remarks in the last section.

\section{Some lemmas}
The Fermat quotient of an integer $a$ with respect to an odd prime $p$ is given by
\begin{align*}
q_p(a)=\frac{a^{p-1}-1}{p},
\end{align*}
which plays an important role in the study of cyclotomic fields.
\begin{lem} \label{lem1}(Eisenstein)
Suppose $p$ is an odd prime and $r$ is a positive integer. For non-zero $p$-adic integers $a$ and $b$, we have
\begin{align*}
q_p(ab)&\equiv q_p(a)+q_p(b)\pmod{p},\\
q_p(a^r)&\equiv rq_p(a) \pmod{p}.
\end{align*}
\end{lem}

\begin{lem}\label{lem2}
(Lehmer \cite{Lehm}) Let $H_n=\sum_{k=1}^n\frac{1}{k}$ be the $n$-th harmonic number. For any prime $p\ge 5$, we have
\begin{align*}
&H_{\lfloor p/2\rfloor}\equiv -2q_p(2)\pmod{p},\quad
H_{\lfloor p/3\rfloor}\equiv -\frac{3}{2}q_p(3)\pmod{p},\\
&H_{\lfloor p/4\rfloor}\equiv -3q_p(2)\pmod{p},\quad
H_{\lfloor p/6\rfloor}\equiv -2q_p(2)-\frac{3}{2}q_p(3)\pmod{p},
\end{align*}
where $\lfloor x \rfloor$ denotes the greatest integer less than or equal to a real number $x$.
\end{lem}

\begin{lem}
If $n$ is a positive integer, then
\begin{align}
&\sum_{k=0}^{n}(-1)^k{n\choose k}{n+k\choose k}=(-1)^n,\label{bb1}\\
&\sum_{k=1}^{n}(-1)^k{n\choose k}{n+k\choose k}H_k=2(-1)^nH_n,\label{bb2}\\
&\sum_{k=1}^{n}(-1)^k{n\choose k}{n+k\choose k}\sum_{i=1}^k\frac{1}{n+i}=(-1)^nH_n.\label{bb3}
\end{align}
\end{lem}
{\it Proof.}
Prodinger \cite{Prod} has given a proof of \eqref{bb1}-\eqref{bb2} by  partial fraction
decomposition and creative telescoping, see also \cite{OS}. Using the same method, Prodinger \cite{Prod} also obtained the following identity:
\begin{align}
\sum_{k=1}^{n}(-1)^k{n\choose k}{n+k\choose k}H_{n+k}=2(-1)^nH_n.\label{bb4}
\end{align}
It follows that
\begin{align*}
&\sum_{k=1}^{n}(-1)^k{n\choose k}{n+k\choose k}\sum_{i=1}^k\frac{1}{n+i}\\
&=\sum_{k=1}^{n}(-1)^k{n\choose k}{n+k\choose k}(H_{n+k}-H_n)\\
&=(-1)^nH_n. \quad \text{(by \eqref{bb1} and \eqref{bb4})}
\end{align*}
This concludes the proof of \eqref{bb3}.
\qed

\begin{lem}\label{lem4}
Let $p\ge 5$ be a prime, $r$ and $0\le k \le p-1$  be non-negative integers. For $x\in \{\frac{1}{2},\frac{1}{3},\frac{1}{4},\frac{1}{6}\}$, we have
\begin{align}
\frac{(x)_{k+rp}(1-x)_{k+rp}}{(1)_{k+rp}^2}&\equiv \frac{(x)_{r}(1-x)_{r}}{(1)_{r}^2}\cdot\frac{(x)_{k}(1-x)_{k}}{(1)_{k}^2}\notag\\
&\times \left(1+2rpH_{\lfloor px \rfloor}-2rpH_k+rp\sum_{i=0}^{k-1}\left(\frac{1}{x+i}+\frac{1}{1-x+i}\right)\right) \pmod{p^2}.\label{bb5}
\end{align}
\end{lem}
{\it Proof.}
Note that
\begin{align*}
\frac{(x)_n(1-x)_n}{(1)_n^2}=\frac{S_x(n)}{a_x^n},
\end{align*}
where
\begin{align*}
S_{x}(n)=\begin{cases}
{2n\choose n}^2,&\text{if $x=\frac{1}{2}$,}\\[5pt]
{2n\choose n}{3n\choose n},&\text{if $x=\frac{1}{3}$,}\\[5pt]
{2n\choose n}{4n\choose 2n},&\text{if $x=\frac{1}{4}$,}\\[5pt]
{3n\choose n}{6n\choose 3n},&\text{if $x=\frac{1}{6}$,}
\end{cases}
\end{align*}
and $a_{1/2}=16, a_{1/3}=27,a_{1/4}=64$ and $a_{1/6}=432$.
We see that \eqref{bb5} is equivalent to
\begin{align}
S_{x}(k+rp)&\equiv a_x^{r(p-1)}S_x(r)S_x(k)\notag\\
&\times\left(1+2rpH_{\lfloor px \rfloor}-2rpH_k+rp\sum_{i=0}^{k-1}\left(\frac{1}{x+i}+\frac{1}{1-x+i}\right)\right) \pmod{p^2}.\label{bb6}
\end{align}

By Lemma \ref{lem1} and \ref{lem2}, we have
\begin{align*}
a_x^{r(p-1)}&\equiv 1+pq_p(a_x^r)\notag\\
&\equiv 1+rpq_p(a_x)\notag\\
&\equiv 1-2rpH_{\lfloor px \rfloor} \pmod{p^2}.
\end{align*}
So in order to prove \eqref{bb6}, it suffices to show that
\begin{align}
S_{x}(k+rp)&\equiv S_x(r)S_x(k)\notag\\
&\times\left(1-2rpH_k+rp\sum_{i=0}^{k-1}\left(\frac{1}{x+i}+\frac{1}{1-x+i}\right)\right) \pmod{p^2}.\label{bb7}
\end{align}

Noting that
\begin{align*}
&\sum_{j=0}^{k-1}\left(\frac{1}{j+1/2}+\frac{1}{j+1/2}\right)=4H_{2k}-2H_k,\\
&\sum_{j=0}^{k-1}\left(\frac{1}{j+1/3}+\frac{1}{j+2/3}\right)=3H_{3k}-H_k,\\
&\sum_{j=0}^{k-1}\left(\frac{1}{j+1/4}+\frac{1}{j+3/4}\right)=4H_{4k}-2H_{2k},\\
&\sum_{j=0}^{k-1}\left(\frac{1}{j+1/6}+\frac{1}{j+5/6}\right)=6H_{6k}-3H_{3k}-2H_{2k}+H_k,
\end{align*}
\eqref{bb7} becomes the following congruences modulus $p^2$
\begin{align}
&{2rp+2k\choose rp+k}^2\equiv {2r\choose r}^2{2k\choose k}^2(1+rp(4H_{2k}-4H_k)),\label{bb8}\\
&{2rp+2k\choose rp+k}{3rp+3k\choose rp+k}
\equiv {2r\choose r}{3r\choose r}{2k\choose k}{3k\choose k}(1+rp(3H_{3k}-3H_k)),\label{bb9}\\
&{2rp+2k\choose rp+k}{4rp+4k\choose 2rp+2k}
\equiv {2r\choose r}{4r\choose 2r}{2k\choose k}{4k\choose 2k}
(1+rp(4H_{4k}-2H_{2k}-2H_k)),\label{bb10}\\
&{3rp+3k\choose rp+k}{6rp+6k\choose 3rp+3k}
\equiv {3r\choose r}{6r\choose 3r}{3k\choose k}{6k\choose 3k}(1+rp(6H_{6k}-3H_{3k}-2H_{2k}-H_k)).\label{bb11}
\end{align}
Next we only prove \eqref{bb8}. The proof of \eqref{bb9}-\eqref{bb11} runs analogously.

Note that
\begin{align}
{2rp+2k\choose rp+k}&={2rp\choose rp}\prod_{i=1}^{2k}(2rp+i)/\prod_{i=1}^{k}(rp+i)^2\notag\\
&\equiv {2r\choose r}\prod_{i=1}^{2k}(2rp+i)/\prod_{i=1}^{k}(rp+i)^2 \pmod{p^2},\label{bb12}
\end{align}
where we have utilized the Babbage's theorem \cite{Bab}
\begin{align*}
{ap\choose bp}\equiv {a\choose b} \pmod{p^2}.
\end{align*}

Now we consider the following rational function:
\begin{align}
f(x)=\prod_{i=1}^{2k}(2rx+i)/\prod_{i=1}^{k}(rx+i)^2.\label{bb13}
\end{align}
Taking the logarithmic derivative on both sides of \eqref{bb13} gives
\begin{align}
\frac{f'(x)}{f(x)}=\sum_{i=1}^{2k}\frac{2r}{2rx+i}-2\sum_{i=1}^{k}\frac{r}{rx+i}.
\label{bb14}
\end{align}
From \eqref{bb13} and \eqref{bb14}, we have
\begin{align*}
f(0)={2k\choose k}\quad \text{and}\quad f'(0)=r{2k\choose k}(2H_{2k}-2H_{k}).
\end{align*}

Now we get the first two terms of the Taylor expansion for $f(x)$:
\begin{align}
f(x)={2k\choose k}+rx{2k\choose k}(2H_{2k}-2H_{k})+\mathcal{O}\left(x^2\right).\label{bb15}
\end{align}
Combining \eqref{bb12} and \eqref{bb15}, we immediately get
\begin{align*}
{2rp+2k\choose rp+k}\equiv {2r\choose r}{2k\choose k}(1+rp(2H_{2k}-2H_k))\pmod{p^2}.
\end{align*}
It follows that
\begin{align*}
{2rp+2k\choose rp+k}^2\equiv {2r\choose r}^2{2k\choose k}^2(1+rp(4H_{2k}-4H_k))\pmod{p^2}.
\end{align*}
This concludes the proof of \eqref{bb8}.
\qed

\begin{lem}\label{lem5}
Suppose $p\ge 5$ is a prime and $0\le k \le p-1$ is an integer. For $x\in\{\frac{1}{2},\frac{1}{3},\frac{1}{4},\frac{1}{6}\}$, we have
\begin{align}
\frac{(x)_k(1-x)_k}{(1)_k^2}\equiv (-1)^k{\lfloor px \rfloor\choose k}{\lfloor px \rfloor+k\choose k}\pmod{p}.\label{cc5}
\end{align}
\end{lem}
{\it Proof.}
It suffices to show that for $0\le k\le p-1$
\begin{align}
(\lfloor px \rfloor+1)_k(-\lfloor px \rfloor)_k
\equiv (x)_k(1-x)_k \pmod{p}.\label{bb17}
\end{align}
For any prime $p\ge 5$, there exists $\varepsilon\in \{1,-1\}$ such that $p\equiv \varepsilon \pmod{2,3,4,6}$. Next we only prove \eqref{bb17} for $x=\frac{1}{3}$. The proof of other three cases runs similarly.

If $p\equiv 1 \pmod{3}$, then $\lfloor p/3 \rfloor=\frac{p-1}{3}$, and hence
\begin{align*}
\left(\frac{p+2}{3}\right)_k\left(\frac{-p+1}{3}\right)_k\equiv \left(\frac{2}{3}\right)_k\left(\frac{1}{3}\right)_k\pmod{p}.
\end{align*}

If $p\equiv -1 \pmod{3}$, then $\lfloor p/3 \rfloor=\frac{p-2}{3}$, and so
\begin{align*}
\left(\frac{p+1}{3}\right)_k\left(\frac{-p+2}{3}\right)_k\equiv \left(\frac{2}{3}\right)_k\left(\frac{1}{3}\right)_k\pmod{p}.
\end{align*}
This implies that \eqref{bb17} is true for $x=\frac{1}{3}$.
\qed

\section{A new proof of Theorem \ref{thmsun}}
Letting $x\to -x$ in Theorem \ref{thmsun}, \eqref{sun-1} is equivalent to
\begin{align}
\pFq{2}{1}{-x,1+x}{1}{1}_p\equiv (-1)^{\langle x\rangle_p} \pmod{p^2}.\label{sun-2}
\end{align}
It is easy to see that
\begin{align*}
\frac{(-x)_k(1+x)_k}{(1)_k^2}=(-1)^k{x\choose k}{x+k\choose k}.
\end{align*}
Let $\delta$ denote the number $\delta=(x-\langle x\rangle_p)/p$. It is clear that $\delta$ is a $p$-adic integer and $x=\langle x\rangle_p+\delta p$. Note that
\begin{align*}
{x\choose k}{x+k\choose k}&={\langle x\rangle_p+\delta p\choose k}{\langle x\rangle_p+\delta p+k\choose k}\notag\\
&=\prod_{i=1}^k(\langle x\rangle_p+\delta p+1-i)\prod_{i=1}^k(\langle x\rangle_p+\delta p+i)\left(\prod_{i=1}^{k}i^{-1}\right)^2\notag\\
&\equiv
{\langle x\rangle_p\choose k}{\langle x\rangle_p+k\choose k}
\left(1+\delta p\left(\sum_{i=1}^k\frac{1}{\langle x\rangle_p+i}+\sum_{i=1}^k\frac{1 }{\langle x\rangle_p+1-i}\right)\right)\pmod{p^2}.
\end{align*}
It follows that
\begin{align}
\text{LHS \eqref{sun-2}}
&\equiv \sum_{k=0}^{p-1}(-1)^k{\langle x\rangle_p\choose k}{\langle x\rangle_p+k\choose k}\notag\\
&\times\left(1+\delta p\left(\sum_{i=1}^k\frac{1}{\langle x\rangle_p+i}+\sum_{i=1}^k\frac{1 }{\langle x\rangle_p+1-i}\right)\right)\pmod{p^2}.\label{sun-3}
\end{align}

Let $b=p-\langle x\rangle_p$. It is clear that $\langle x\rangle_p\equiv -b\pmod{p}$ and $0\le b-1\le p-1$. By \eqref{bb3}, we have
\begin{align}
&\sum_{k=0}^{p-1}(-1)^k{\langle x\rangle_p\choose k}{\langle x\rangle_p+k\choose k}\sum_{i=1}^k\frac{1 }{\langle x\rangle_p+1-i}\notag\\
&\equiv -\sum_{k=0}^{p-1}{-b\choose k}{-b+k\choose k}\sum_{i=1}^k\frac{1 }{b-1+i}\pmod{p}\notag\\
&=-\sum_{k=0}^{p-1}{b-1\choose k}{{b-1+k\choose k}}\sum_{i=1}^k\frac{1 }{b-1+i}\notag\\
&=(-1)^bH_{b-1}\notag\\
&\equiv(-1)^{\langle x\rangle_p+1}H_{\langle x\rangle_p}\pmod{p},\label{sun-4}
\end{align}
where we have used the fact that ${-b\choose k}{-b+k\choose k}={b-1\choose k}{{b-1+k\choose k}}$
in the second step and $H_{p-k-1}\equiv H_k\pmod{p}$ for $0\le k\le p-1$ in the last step.

Using \eqref{bb1} and \eqref{bb3}, we obtain
\begin{align}
\sum_{k=0}^{p-1}(-1)^k{\langle x\rangle_p\choose k}{\langle x\rangle_p+k\choose k}=(-1)^{\langle x\rangle_p},\label{sun-5}
\end{align}
and
\begin{align}
\sum_{k=0}^{p-1}(-1)^k{\langle x\rangle_p\choose k}{\langle x\rangle_p+k\choose k}\sum_{i=1}^k\frac{1 }{\langle x\rangle_p+i}=(-1)^{\langle x\rangle_p}H_{\langle x\rangle_p}.\label{sun-6}
\end{align}
Substituting \eqref{sun-4}-\eqref{sun-6} into \eqref{sun-3}, we complete the proof of \eqref{sun-2}.

\section{Proof of Theorem \ref{thm2}}
Assume $x\in\{\frac{1}{2},\frac{1}{3},\frac{1}{4},\frac{1}{6}\}$.
We first prove that for any non-negative integer $r$
\begin{align}
\sum_{k=rp}^{(r+1)p-1}\frac{(x)_k(1-x)_k}{(1)_k^2}\equiv \frac{(x)_{r}(1-x)_{r}}{(1)_{r}^2}\cdot\pFq{2}{1}{x,1-x}{1}{1}_p \pmod{p^2}. \label{cc1}
\end{align}
Letting $k\to k+rp$ on the left-hand side of \eqref{cc1} and then applying Lemma \ref{lem4}, we obtain
\begin{align}
\sum_{k=rp}^{(r+1)p-1}\frac{(x)_k(1-x)_k}{(1)_k^2}
&\equiv\frac{(x)_{r}(1-x)_{r}}{(1)_{r}^2}\sum_{k=0}^{p-1}\frac{(x)_{k}(1-x)_{k}}{(1)_{k}^2}\notag\\
&\times \left(1+2rpH_{\lfloor px \rfloor}-2rpH_k+rp\sum_{i=0}^{k-1}\left(\frac{1}{x+i}+\frac{1}{1-x+i}\right)\right) \pmod{p^2}.\label{cc2}
\end{align}

Using the following identity \cite[Theorem 1]{Taur2}:
\begin{align*}
\frac{(x)_k(1-x)_k}{(1)_k^2}
\sum_{i=0}^{k-1}\left(\frac{1}{x+i}+\frac{1}{1-x+i}\right)=\sum_{i=0}^{k-1}
\frac{(x)_i(1-x)_i}{(1)_i^2}\cdot\frac{1}{k-i},
\end{align*}

we have
\begin{align}
&\sum_{k=0}^{p-1}\frac{(x)_{k}(1-x)_{k}}{(1)_{k}^2}\sum_{i=0}^{k-1}\left(\frac{1}{x+i}+\frac{1}{1-x+i}\right)\notag\\
&=\sum_{k=0}^{p-1}\sum_{i=0}^{k-1}\frac{(x)_i(1-x)_i}{(1)_i^2}\cdot\frac{1}{k-i}\notag\\
&=\sum_{i=0}^{p-2}\frac{(x)_i(1-x)_i}{(1)_i^2}\sum_{k=i+1}^{p-1}\frac{1}{k-i}\notag\\
&\equiv\sum_{i=0}^{p-1}\frac{(x)_i(1-x)_i}{(1)_i^2}H_i \pmod{p},\label{cc3}
\end{align}
since $H_{p-1-i}\equiv H_i \pmod{p}$ and $(x)_i(1-x)_i\equiv 0\pmod{p}$ for $i=p-1$.

It follows from \eqref{cc2} and \eqref{cc3} that
\begin{align}
&\sum_{k=rp}^{(r+1)p-1}\frac{(x)_k(1-x)_k}{(1)_k^2}\notag\\
&\equiv \frac{(x)_{r}(1-x)_{r}}{(1)_{r}^2}\left( \sum_{k=0}^{p-1}\frac{(x)_{k}(1-x)_{k}}{(1)_{k}^2}
+rp\sum_{k=0}^{p-1}\frac{(x)_{k}(1-x)_{k}}{(1)_{k}^2}(2H_{\lfloor px\rfloor}-H_k)\right)\pmod{p^2}.\label{cc4}
\end{align}
By \eqref{cc5}, we have
\begin{align}
&\sum_{k=0}^{p-1}\frac{(x)_{k}(1-x)_{k}}{(1)_{k}^2}(2H_{\lfloor px\rfloor}-H_k)\notag\\
&\equiv \sum_{k=0}^{p-1}(-1)^k{\lfloor px\rfloor\choose k}{\lfloor px\rfloor+k\choose k}(2H_{\lfloor px\rfloor}-H_k)\pmod{p}\notag\\
&=0 .\quad\text{(by \eqref{bb1} and \eqref{bb2})}\label{cc6}
\end{align}
Substituting \eqref{cc6} into \eqref{cc4}, we complete the proof of \eqref{cc1}.

Taking the sum over $r$ from $0$ to $n-1$ on both sides of \eqref{cc1} gives
\begin{align}
\pFq{2}{1}{x,1-x}{1}{1}_{np}&\equiv \pFq{2}{1}{x,1-x}{1}{1}_p\cdot\pFq{2}{1}{x,1-x}{1}{1}_n \pmod{p^2}.\label{cc7}
\end{align}
The proof of Theorem \ref{thm2} then follows from \eqref{sun-1} and \eqref{cc7}.

\section{Concluding remarks}
Numerical calculation suggests that supercongruence \eqref{rv} can not be extended to any $p$-adic integer $x$ in the direction of Theorem \ref{thmsun}.

Recently, Z.-W. Sun \cite[Conjecture 5.4]{Szw2} made four challenging conjectures which extend Theorem \ref{thm2} and some results proved by Z.-H. Sun \cite{Szh2} and Z.-W. Sun \cite{Szw}.

\begin{conj}\label{conj}(Z.-W. Sun)
Let $p\ge 5$ be a prime and $n$ be a positive integer. Then
\begin{align*}
\frac{16^n}{n^2{2n\choose n}^2}\left(\sum_{k=0}^{np-1}\frac{{2k\choose k}^2}{16^k}-\left(\frac{-1}{p}\right)\sum_{k=0}^{n-1}\frac{{2k\choose k}^2}{16^k}\right)
&\equiv -4p^2E_{p-3}\pmod{p^3},\\
\frac{27^n}{n^2{2n\choose n}{3n\choose n}}\left(\sum_{k=0}^{np-1}\frac{{2k\choose k}{3k\choose k}}{27^k}-\left(\frac{-3}{p}\right)\sum_{k=0}^{n-1}\frac{{2k\choose k}{3k\choose k}}{27^k}\right)
&\equiv -\frac{3}{2}p^2B_{p-2}\left(\frac{1}{3}\right)\pmod{p^3},\\
\frac{64^n}{n^2{2n\choose n}{4n\choose 2n}}\left(\sum_{k=0}^{np-1}\frac{{2k\choose k}{4k\choose 2k}}{64^k}-\left(\frac{-2}{p}\right)\sum_{k=0}^{n-1}\frac{{2k\choose k}{4k\choose 2k}}{64^k}\right)
&\equiv -p^2E_{p-3}\left(\frac{1}{4}\right)\pmod{p^3},\\
\frac{432^n}{n^2{3n\choose n}{6n\choose 3n}}\left(\sum_{k=0}^{np-1}\frac{{3k\choose k}{6k\choose 3k}}{432^k}-\left(\frac{-1}{p}\right)\sum_{k=0}^{n-1}\frac{{3k\choose k}{6k\choose 3k}}{432^k}\right)
&\equiv -20p^2E_{p-3}\pmod{p^3},
\end{align*}
where $E_m$ is the $m$-th Euler number, $E_m(x)$ and $B_m(x)$ denote the Euler polynomial and the Bernoulli polynomial of degree $m$, respectively.
\end{conj}

Noting that
\begin{align*}
&\frac{{2k\choose k}^2}{16^k}=\frac{\left(\frac{1}{2}\right)_k^2}{(1)_k^2},\quad\quad
\frac{{2k\choose k}{3k\choose k}}{27^k}=\frac{\left(\frac{1}{3}\right)_k \left(\frac{2}{3}\right)_k}{(1)_k^2},\\
&\frac{{2k\choose k}{4k\choose 2k}}{64^k}=\frac{\left(\frac{1}{4}\right)_k \left(\frac{3}{4}\right)_k}{(1)_k^2},\quad\quad
\frac{{3k\choose k}{6k\choose 3k}}{432^k}=\frac{\left(\frac{1}{6}\right)_k \left(\frac{5}{6}\right)_k}{(1)_k^2},
\end{align*}
we can directly deduce Theorem \ref{thm2} from Conjecture \ref{conj}. Unfortunately, the
method in this paper is not applicable for proving Conjecture \ref{conj}. We hope that the interested reader can work on this problem.

\vskip 5mm \noindent{\bf Acknowledgments.} The author would like to
thank Professor V.J.W. Guo for valuable comments which improved the presentation of the paper.

\end{document}